\documentclass[a4paper, BCOR=0.0mm, DIV=calc]{scrartcl}
\usepackage[T1]{fontenc}
\usepackage[utf8]{inputenc}
\usepackage[ngerman]{babel}
\usepackage[osf,sc]{mathpazo}
\usepackage{textcomp}
\usepackage{microtype}
\usepackage{amsmath, amssymb, amsfonts}
\usepackage{tabularx}
\KOMAoptions{twoside=false, twocolumn=false, headinclude=false, footinclude=false, mpinclude=false, pagesize=auto}
\recalctypearea

\setkomafont{section}{\mdseries\scshape\Large}
\setkomafont{title}{\mdseries\scshape}

\begin{document}
\title{Entwicklung des unendlichen Produkts \\ {\Large $(1-x)(1-xx)(1-x^3)(1-x^4)(1-x^5)(1-x^6)\mathrm{etc}$} \\ in eine einfache Reihe\footnote{
Originaltitel: "`Evolutio producti infiniti $(1-x)(1-xx)(1-x3)(1-x4)(1-x5) etc.$ in seriem simplicem"', erstmals publiziert in "`\textit{Acta Academiae Scientarum Imperialis Petropolitinae} 1780, 1783, pp. 47-55"', Nachdruck in "`\textit{Opera Omnia}: Series 1, Volume 3, pp. 472 - 479"', Eneström-Nummer E541, übersetzt von: Alexander Aycock, Textsatz: Artur Diener,  im Rahmen des Projektes "`Eulerkreis Mainz"'}}
\author{Leonhard Euler}
\date{}
\maketitle
\paragraph{§1}
Nachdem	$s = (1-x)(1-xx)(1-x^3)(1-x^4)\mathrm{etc}$ gesetzt worden ist, ist leicht klar, dass gelten wird:
\[
	s = 1 -x -xx(1-x)-x^3(1-x)(1-xx)-x^4(1-x)(1-xx)(1-x^3) - \mathrm{etc};
\]
weil diese Reihe schon eine unendliche ist, fragt man, wenn ihre einzelnen Terme entwickelt werden, welche Reihe nach den einfachen Potenzen von $x$ hervorgehen wird. Weil also die ersten zwei Terme $1-x$ schon entwickelt sind, schreibe man anstelle aller übrigen den Buchstaben $A$, sodass $s = 1-x-A$ ist, und daher
\[
	A = xx(1-x) + x^3(1-x)(1-xx) + x^4(1-x)(1-xx)(1-x^3) + \mathrm{etc.}
\]
\paragraph{§2}
Weil ja all diese Terme den gemeinsamen Faktor $1-x$ haben, teile man, nachdem er entwickelt worden ist, die einzelnen Terme in zwei Teile auf, welche wir so darstellen wollen:
\begin{alignat*}{4}
A &= xx + &x^3(1-xx)+&x^4(1-xx)(1-x^3) + &x^5(1-xx)(1-x^3)(1-x^4) \\
  &- x^3 - &x^4(1-xx)-&x^5(1-xx)(1-x^3) - &x^6(1-xx)(1-x^3)(1-x^4)
\end{alignat*}
Daher fasse man gleich die beiden Teile, die mit derselben Potenz von $x$ versehen sind, zu einem zusammen, und es wird für $A$ die folgende Form resultieren:
\begin{align*}
A = xx &-x^5 -x^7(1-xx)-x^9(1-xx)(1-x^3) \\
  &- x^{11}(1-xx)(1-x^3)(1-x^4) - \mathrm{etc,}
\end{align*}
wo die ersten beiden Terme $xx-x^5$ schon entwickelt sind; die folgenden schreiten durch diese Potenzen fort: $x^7,\, x^9,\, x^{11},\, x^{13},\, x^{15}$, deren Exponenten um zwei wachsen.
\paragraph{§3}
Wir wollen nun auf dieselbe Weise wie zuvor	$A = xx - x^5 - B$ setzen, sodass
\begin{align*}
B = +x^7&(1-xx) + x^9(1-xx)(1-x^3) \\
  + x^{11}&(1-xx)(1-x^3)(1-x^4) + \mathrm{etc.}
\end{align*}
ist, alle Terme welche den gemeinsamen Faktor $1-xx$ haben, nach Entwicklung welches man die einzelnen Terme wie folgt in zwei Teile aufteile:
\begin{alignat*}{4}
 B &= x^7 + x^9&(1-x^3) + &x^{11}(1-x^3)(1-x^4) + x^{13}&(1-x^3)(1-x^4)(1-x^5) \\
   &- x^9 - x^{11}&(1-x^3) - &x^{13}(1-x^3)(1-x^4) - x^{15}&(1-x^3)(1-x^4)(1-x^5). 
\end{alignat*}
Hier fasse man wiederum zwei Terme, die dieselbe Potenz von $x$ vorangestellt haben, zu einem zusammen und es wird hervorgehen:
\begin{align*}
B = x^7 &- x^{12} - x^{15}(1-x^3) - x^{18}(1-x^3)(1-x^4) \\
  &- x^{21}(1-x^3)(1-x^4)(1-x^5) - \mathrm{etc,}
\end{align*}
wo die Potenzen von $x$ schon um $3$ wachsen.
\paragraph{§4}
Man setze nun weiter $B = x^7 - x^{12} - C$, sodass
\begin{align*}
C = x^{15}&(1-x^3)+x^{18}(1-x^3)(1-x^4) \\
  + x^{21}&(1-x^3)(1-x^4)(1-x^5) + \mathrm{etc}
\end{align*}
ist und man löse gleich die einzelnen Terme durch Entwicklung des Faktors $1-x^3$ in zwei Teile auf, und es wird
\begin{alignat*}{4}
C &= x^{15} + &x^{18}(1-x^4) + &x^{21}(1-x^4)(1-x^5) + &x^{24}(1-x^4)(1-x^5)(1-x^6) \\
  &- x^{18} - &x^{21}(1-x^4) - &x^{24}(1-x^4)(1-x^5) - &x^{27}(1-x^4)(1-x^5)(1-x^6)
\end{alignat*}
werden, wo erneut die Glieder, denen die dieselbe Potenz von $x$ vorangestellt ist, zu einem zusammengefasst
\begin{align*}
C = x^{15} - x^{22} &- x^{26}(1-x^4) - x^{30}(1-x^4)(1-x^5) \\
&- x^{24}(1-x^4)(1-x^5)(1-x^6)\mathrm{etc}
\end{align*}
liefern werden, wo die vorangestellten Potenzen um $4$ wachsen.
\paragraph{§5}
Man setze nun $C = x^{15} - x^{22} - D$, sodass
\begin{align*}
	D = x^{26}(1-x^4) + &x^{30}(1-x^4)(1-x^5) \\
					  +	&x^{34}(1-x^4)(1-x^5)(1-x^6)\mathrm{etc} \\
\end{align*}
ist, welche Terme man durch die Entwicklung des Faktors $1-x^4$ auf diese Weise in zwei Teile aufteile
\begin{alignat*}{4}
D = &x^{26} + &x^{30}(1-x^5) + &x^{34}(1-x^5)(1-x^6) + x^{38}&(1-x^5)(1-x^6)(1-x^7) \\
   -&x^{30} - &x^{34}(1-x^5) - &x^{38}(1-x^5)(1-x^6) - x^{42}&(1-x^5)(1-x^6)(1-x^7). 
\end{alignat*}
Indem man die zwei nun wie bisher aufteilt, berechnet man
\begin{align*}
D = x^{26} - x^{35} - &x^{40}(1-x^5) - x^{45}(1-x^5)(1-x^6) \\
                    - &x^{50}(1-x^5)(1-x^6)(1-x^7)\mathrm{etc.}
\end{align*}
Hier wachsen also die Potenzen von $x$ um $5$.
\paragraph{§6}
Man setze $D = x^{26} - x^{35} - E$, sodass
\begin{align*}
E = x^{40}(1-x^5) + &x^{45}(1-x^5)(1-x^6) \\
                  + &x^{50}(1-x^5)(1-x^6)(1-x^7)\mathrm{etc}
\end{align*}
ist und es geht, nachdem die Auflösung in zwei Teile wie bisher ausgeführt wurde,
\begin{alignat*}{3}
E = x^{40} + &x^{45}(1-x^6) + &x^{50}(1-x^6)(1-x^7) + &x^{55}(1-x^6)(1-x^7)(1-x^8) \\
   -x^{45} - &x^{50}(1-x^6) - &x^{55}(1-x^6)(1-x^7) - &x^{60}(1-x^6)(1-x^7)(1-x^8)
\end{alignat*}
hervor. Nachdem aber diese zwei Terme zu einem zusammengefasst worden sind, wird
\begin{align*}
E = x^{45} - x^{51} - &x^{57}(1-x^6) - x^{63}(1-x^6)(1-x^7) \\
                    - &x^{69}(1-x^6)(1-x^7)(1-x^8)\mathrm{etc}
\end{align*}
sein, wo die Potenzen von $x$ um $6$ wachsen.
\paragraph{§7}
Weil das Bidlungsgesetz, nach welchem diese Operationen weiter fortzusetzen sind, hinreichend klar ist, werden wir, wenn die letzten für die einzelnen Buchstaben $A$, $B$, $C$, $D$ gefundenen Werte der Reihe nach eingesetzt werden, für die gesuchte Reihe die folgende Form finden:
\begin{align*}
	S =& 1 - x - xx + x^5 + x^7 - x^{12} - x^{15} + x^{20} \\
	  &+ x^{26} - x^{35} - x^{40} + x^{51} + \mathrm{etc}.
\end{align*}
Hier geht nun also die ganze Frage darauf zurück, dass die Ordnung bestimmt wird, nach welcher die Exponenten der Potenzen von $x$ immer weiter vermehrt werden, weil ja aus den durchgeführten Rechnungen schon hinreichend klar ist, dass die Zeichen "`$+$"' und "`$-$"' abwechselnd so folgen, dass sie je doppelt auftreten.
\paragraph{§8}
Um also dieses Bildungsgesetz leichter zu untersuchen, wollen wir sehen, wie in diesen Buchstaben diese Zahlen entstanden sind. Für dieses Ziel wollen wir zumindest die ersten Terme jedes Buchstabens, welche in ihrer ersten Form ausgedrückt worden sind, der Reihe nach anordnen:
\begin{alignat*}{4}
A &= xx(1-x) \qquad \quad 7 &= 3 + 4 &= 4 + 1 + 3 &= 3 +1 +1 +2 \\
B &= x^7(1-xx) \qquad 15 &= 4 + 11 &= 4 + 2  + 9 &= 4 +2 +2 +7 \\
C &= x^{15}(1-x^3) \qquad 26 &= 5 + 21 &= 5 + 3 + 18 &= 5 + 3 + 3 + 15 \\
D &= x^{26}(1-x^4) \qquad 40 &= 6 + 34 &= 6 + 4 + 30 &= 6 + 4 + 4 + 26 \\
E &= x^{40}(1-x^5) \qquad 57 &= 7 + 50 &= 7 + 5 + 45 &= 7 + 5 + 5 + 40 \\
&\mathrm{etc} &  \mathrm{etc.} & &
\end{alignat*}
Hier sehen wir natürlich aus der Entwicklung des Buchstaben $A$, dass die Zahl $7$ aus dem Aggregat $3+4$ entsteht, dann in der Tat $4$ aus $1+3$ entsteht und schließlich $3$ aus $1+2$, welche Auflösung also
\[
	7 = 3 + 4 = 3 + 1 + 3 = 3+1+1+2
\]
geben wird. Und dieselbe Struktur ist in den folgenden Buchstaben beobachtet worden, wo die letzten Zahlen nach der Ordnung $2$, $7$, $15$, $26$, $40$ fortschreiten.
\paragraph{§9}
Daraus ist schon klar, dass die Differenzen der Zahlen $2$, $7$, $15$, $40$, $57$, etc eine arithemtische Progression festsetzen, woher der allgemeine Term dieser Zahlen
\[
	2 + 5(n-1) + \frac{2(n-1)(n-2)}{1\cdot 2} = \frac{3nn+n}{2}
\]
sein wird. Die Exponenten aber, die diesen vorausgehen, waren $1$, $5$, $12$, $22$, $39$, $51$ von den Zahlen $1$, $2$, $3$, $4$, $5$ und im Allgemeinen von der Zahl $n$ selbst, sodass der Exponent, der der Formel $\frac{3nn+n}{2}$ vorausgeht, $\frac{3nn-n}{2}$ sein wird.
\paragraph{§10}
Nun erkennen wir also vollkommen die einfache gefundene Reihe, die dem vorgelegten unendlichen Produkt
\[
	(1-x)(1-xx)(1-x^3)(1-x^4)\mathrm{etc}
\]
gleich ist. Weil nämlich diese Reihe gefunden worden ist,
\begin{align*}
	s = 1 - x &- x^2 + x^5 + x^7 - x^{12} - x^{15} + x^{22} + x^{26} \\
	  &- x^{35} - x^{40} + x^{51} + \mathrm{etc},
\end{align*}
sind wir nun sicher, dass in ihr andere Potenzen von $x$ nicht auftauchen, wenn deren Exponenten nicht in dieser Form, $\frac{3nn\pm n}{2}$, enthalten sind, und zwar so, dass wenn $n$ eine ungerade Zahl war, die beiden daher entstehenden Terme das Vorzeichen "`$-$"' haben werden, die aber aus den geraden entstehen, das Vorzeichen "`$+$"'.
\section*{Andere Untersuchung derselben Reihe}
\paragraph{§11}
Dieselbe Reihe, die nach den Potenzen von $x$ fortschreitet, kann auch auf folgende Weise untersucht werden. Weil natürlich
\begin{align*}
	s = 1 - x - xx&(1-x) - x^3(1-x)(1-xx) \\
	          -x^4&(1-x)(1-xx)(1-x^3)\mathrm{etc}
\end{align*}
ist, entwickle man sofort das zweite Glied $-xx(1-x)$, sodass
\begin{align*}
	s = 1 - x &-xx + x^3 - x^3(1-x)(1-xx) \\
	  &-x^4(1-x)(1-xx)(1-x^3)\mathrm{etc}
\end{align*}
ist und setze 
\[
	s = 1 - x - xx + A,
\]
dass
\begin{align*}
	A = x^3 - x^3&(1-x)(1-xx) \\
			- x^4&(1-x)(1-xx)(1-x^3)-\mathrm{etc}
\end{align*}
ist, deren einzelne Glieder man durch die Entwicklung des Faktors $1-x$ in zwei Teile aufteile, dass
\begin{alignat*}{3}
	A = x^3 - x^3&(1-xx) - x^4&(1-xx)(1-x^3) - x^5&(1-xx)(1-x^3)(1-x^4) \\
	        + x^4&(1-xx) + x^5&(1-xx)(1-x^3) + x^6&(1-xx)(1-x^3)(1-x^4)
\end{alignat*}
hervorgeht. Hier werden wiederum die mit der selben Potenz von $x$ versehenen Glieder zusammengefasst
\begin{align*}
	A = x^5 + x^7&(1-xx) + x^9(1-xx)(1-x^3) \\
	       x^{11}&(1-xx)(1-x^3)(1-x^4)\mathrm{etc}
\end{align*}
liefern.
\paragraph{§12}
Hier entwickle man das zweite Glied, dass
\begin{align*}
	A = x^5 &+ x^7 - x^9 + x^9(1-xx)(1-x^3) \\
	&+x^{11}(1-x^2)(1-x^3)(1-x^4)\mathrm{etc}
\end{align*}
hervorgeht. Man setze gleich $A = x^5 + x^7 - B$, dass
\[
B = x^9 - x^9(1-xx)(1-x^3) - x^{11}(1-xx)(1-x^3)(1-x^4)\mathrm{etc}
\]
ist; wenn daher überall der Faktor $1-xx$ entwickelt wird, wird man
\begin{alignat*}{3}
	B = x^5 - x^9&(1-x^3) - x^{11}&(1-x^3)(1-x^4) - x^{13}&(1-x^3)(1-x^4)(1-x^5) \\
            + x^{11}&(1-x^3) + x^{13}&(1-x^3)(1-x^4) + x^{15}&(1-x^3)(1-x^4)(1-x^5)
\end{alignat*}
erhalten, dann wird aber durch Zusammenfassen der ersten beiden Glieder
\begin{align*}
	B = x^{12} + x^{15}&(1-x^3) + x^{18}(1-x^3)(1-x^4) \\
	           + x^{21}&(1-x^3)(1-x^4)(1-x^5)\mathrm{etc}
\end{align*}
entstehen.
\paragraph{§13}
Man entwickle in gleicher Weise das zweite Glied und setze $B = x^{12} + x^{15} - C$; und es wird
\[
	C  = x^{18} - x^{18}(1-x^3)(1-x^4) - x^{21}(1-x^3)(1-x^4)(1-x^5)-\mathrm{etc}
\]
sein. Nun entwickle man den zweiten Faktor $1-x^3$, und es wird
\begin{alignat*}{3}
	C = x^{18} - x^{18}&(1-x^4) - x^{21}&(1-x^4)(1-x^5) - x^{24}&(1-x^4)(1-x^5)(1-x^6) \\
	           + x^{21}&(1-x^4) + x^{24}&(1-x^4)(1-x^5) + x^{27}&(1-x^4)(1-x^5)(1-x^6)
\end{alignat*}
werden. Daher wird durch Zusammenfassen zweier Glieder
\[
	C = x^{22} + x^{26}(1-x^4) + x^{30}(1-x^4)(1-x^5) + x^{34}(1-x^4)(1-x^5)(1-x^6) + \mathrm{etc}
\]
werden.
\paragraph{§14}
Nachdem nun hier wiederum das zweite Glied entwickelt worden ist, setze man $C = x^{22}+x^{26} - D$, und es wird
\[
	D = x^{30} - x^{30}(1-x^4)(1-x^5) - x^{34}(1-x^4)(1-x^5)(1-x^6)\mathrm{etc}
\]
sein, wo die Entwicklung des Faktors
\begin{alignat*}{3}
	D = x^{30} - x^{30}&(1-x^5) - x^{34}&(1-x^5)(1-x^6) - x^{38}&(1-x^5)(1-x^6)(1-x^7) \\
	           + x^{34}&(1-x^5) + x^{38}&(1-x^5)(1-x^6) + x^{42}&(1-x^5)(1-x^6)(1-x^7)
\end{alignat*}
ergeben wird. Daher wird durch Zusammenfassen zweier Glieder
\[
	D = x^{35} + x^{40}(1-x^5) + x^{45}(1-x^5)(1-x^6) + x^{50}(1-x^5)(1-x^6)(1-x^7)+\mathrm{etc}
\]
werden.
\paragraph{§15}
Nachdem das zweite Glied entwickelt worden ist, setze man erneut $D = x^{35} + x^{40} - E$, und es wird 
\[
	E = x^{45} - x^{45}(1-x^5)(1-x^6) - x^{50}(1-x^5)(1-x^6)(1-x^7)
\]
sein und es wird, nach der Entwicklung des zweiten Faktors $1-x^5$
\begin{alignat*}{3}
E = x^{45} - x^{45}&(1-x^6) - x^{50}&(1-x^6)(1-x^7) - x^{55}&(1-x^6)(1-x^7)(1-x^8) \\
           + x^{50}&(1-x^6) + x^{55}&(1-x^6)(1-x^7) + x^{60}&(1-x^6)(1-x^7)(1-x^8)
\end{alignat*}
werden und durch Zusammenfassen zweier Terme findet man
\begin{align*}
	E = x^{51} + x^{57}&(1-x^6) + x^{63}(1-x^6)(1-x^7) \\
			     x^{69}&(1-x^6)(1-x^7)(1-x^8)\mathrm{etc}
\end{align*}
\paragraph{§16}
Nachdem also diese Werte der Buchstaben $A$, $B$, $C$, $D$, $E$, $F$ gefunden worden sind, wird, wenn die einzelnen nacheinander eingesetzt werden, diese Reihe resultieren:
\[
	1-x-xx+x^5+x^7-x^{12}-x^{15}+x^{22}+x^{26}-x^{35}-x^{40}+\mathrm{etc.}
\]
Hier durchschaut man aber die Struktur der Exponenten leichter: Weil nämlich in den zuerst der Buchstaben $A$, $B$, $C$, $D$ festgesetzten Werten die ersten einfachen Terme $x^3$, $x^9$, $x^{18}$, $x^{30}$, $x^{45}$, etc waren, sind die Exponenten natürlich das Deifache der Dreieckszahlen, woher allgemein für die Zahlen $n$ dieser Exponent $\frac{3nn+n}{2}$ sein wird. Aber diese Terme folgen zwei vorhergehenden Potenzen von $x$ durch dieselbe Differenz $n$, woher die durch zweimaliges Abziehen dieser Zahl von dieser Formel werden die zwei in die gesuchte Reihe eingehenden Potenz entstehen, deren Exponenten als logische Konsequenz
\[
	\frac{3nn+n}{2} \quad \text{und} \quad \frac{3nn-n}{2}
\]
sein werden.
\paragraph{§17}
Daher ist also andererseits klar, dass die ins Unendliche festgesetzte Reihe
\[
	1-x-xx+x^5+x^7-x^{12}-x^{15}+x^{22}+x^{26}-x^{35}-\mathrm{etc}
\]
unendlich Faktoren haben wird, die natürlich $1-x$, $1-xx$, $1-x^3$, $1-x^4$, $1-x^5$, etc sein werden, sodass, wenn zuerst durch $1-x$ geteilt wird, dann in der Tat der Quotient durch $1-xx$, dieser Quotient weiter durch $1-x^3$, und wenn auf diese Weise die Division ins Unendliche fortgesetzt wird, der letzte resultierende Quotient der Einheit gleich sein muss.
\paragraph{§18}
Wenn also diese ins Unendliche laufende Gleichung vorgelegt war:
\[
	1-x-xx+x^5+x^7-x^{12}-x^{15}+x^{22}+x^{26}-\mathrm{etc} = 0,
\]
können all ihre Wurzeln leicht angegeben werden. Zuerst wird nämlich $x=1$ eine Wurzel sein, darauf die beiden Quadratwurzeln aus der Einheit, dann in der Tat die drei Kubikwurzeln aus der Einheit, weiter die vier Biquadratwurzeln aus der Einheit; und auf ähnliche Weise die fünf Quintikwurzeln aus der Einheit, und so weiter, unter welchen also die Einheit selbst unendlich auftaucht; aber in der Tat wird man $-1$ dort finden, wo die Wurzel einer geraden Potenz zu ziehen ist.
\end{document}